\documentclass[a4paper,10pt]{amsart}
\usepackage{amssymb,epsfig,amsmath,amsthm}
\usepackage[font=small,format=plain,labelfont=bf,up,textfont=it,up]{caption}
\usepackage{graphicx,subfig}
\usepackage{float}

\graphicspath{{./Figs/}}
\title{On joint distribution of range and terminal value of a Brownian motion}
\author{Oleg Svirschi}

\begin{document}

\begin{abstract}
In this note, we present the closed form solution for the joint distribution of the range and terminal value of a Brownian motion. Based on this distribution we build a range scaled terminal value distribution and show the derivation steps of its density, further s-density. Finally, we sample the s-density from different groups of currency pairs and compare them with theoretical result.
\end{abstract}
\maketitle

\newtheorem{theorem}{Theorem}[section]
\newtheorem{la}[theorem]{Lemma}
\newtheorem{cor}[theorem]{Corollary}
\newtheorem{remark}[theorem]{Remark}
\newtheorem{prob}[theorem]{Problem}	
\newtheorem{definition}[theorem]{Definition}
\newtheorem{alg}[theorem]{Algorithm}
\newtheorem{prop}[theorem]{Proposition}
\newtheorem{assumption}[theorem]{Assumption}
\newtheorem{example}[theorem]{Example}

\section{Introduction}

The joint distribution of minimum, maximum and terminal value of a Brownian motion is well-know concept, see \cite{trivariate}. The slight confusion comes with the fact that it is presented in a form that it is not a density nor a cumulative density. To get to the density, which is a very useful to have, one needs to derivate that expression further. In this paper, first, we obtain the closed form for the trivariate density of minimum, maximum and terminal value. From this trivariate density, we derive the joint density of terminal value and the range, the difference between maximum and minimum. Next, we obtain the density of the terminal value scaled by the range (s-density to keep it short):

 \[s-stat=\frac{W_T}{W_{max}-W_{min}}\]
\par
\hfill 
\par
 and finally compare the obtained density to the density obtained empirically using FX market data.
\setlength{\parindent}{8ex}

\par

\section{Motivation}
The knowledge about maximum and minimum (or high and low) of a price of an asset is an important information. In a directly traded instrument, extremums are given purely by the fact of transacting at these prices. However in a synthetic instrument, which prices are derived from directly traded prices, the high and the lows are estimated as the instrument is not directly traded. In this work we propose a statistical test for high and lows under different markets. And show some examples from synthetic FX crosses.

\subsection{Structural test for market data quality} When we expect the market data to have some sort of approximate behaviour to a Brownian motion, the s-stat for the data will look very similar to the s-stats provided by the close form solution or empirical simulation. However, if the data has a structural problems, for examples highs are overestimated and lows are underestimated, then such deviations can be shown on the s-stat. In other words, s-stat can help testing the minimum and maximum whether they are truly belong to the process.  Such problems are abundant in the synthetic market data, where high and lows are not directly observable. To calculate the high and the low of a spread or more complicated structure, the data provider needs all possible observations for each of the component of the synthetic instrument. Even taken with the most high frequency, the highs will tend to be underestimated and the lows will tend to be overestimated.  

\subsection{Machine learning} In machine learning, normalization of features by scaling them using previous range is a standard technique which helps training neural networks faster and as well as preventing them from overfitting. Inspecting the scaled features using s-stat can help getting insights into the data a machine learning algorithms being trained with. 

\subsection{Trend and Mean Reversion testing} Comparison of two s-stats can lead to hypothesis testing regarding trend vs mean reversion, while trending series will have more weights on the tails of s-density, mean reverting series will have more weight in the middle.  Structural test for integrated process - separating pure noise from signal has been a problem for financial market professionals for a long time. Structurally, noise (while noise) should have s-stat which corresponds to s-stat of a Brownian motion if taken in cumulative. Deviations from this process can be picked up by s-stat. Most of distributions related to financial markets are unbounded, this makes it difficult to deal with fat tails, by range scaling we make the distribution density bounded, this allows more tangible measurement as well as treatment of outliers.

\section{Derivation of the densities}

\subsection{Density of the minimum, maximum, and terminal value for a standard Brownian motion} For a Brownian motion ${W_t, t>0}$, $W_0=0$: the density of the minimum $l_T=min\{W_t\}_{t=0}^T$, maximum $h_T=max\{W_t\}_{t=0}^T$ and the terminal value $W_T$, can be presented as:

\[P(-l \le l_T \le h_T \le h; W_T\in [x,x+dx])=\psi_{h,l,T}(x)dx\]
\[\psi_{h,l,T}(x)=1_{-l\le x\le h}\sum^{+\infty }_{n=-\infty }{(g_{T}(2n(h+l)-x)-g_{T}(2n(h+l)+x-2h))},\]

\[g_T(x)=\frac{1}{\sqrt{2\pi T  }}e^{-\frac{x^2}{2 T}}\]

\noindent Where we use \cite{feller} p.430, to get to the density from this conditional formulation we derivate in respect to ${l}$ and ${h}$.

\[\psi''(T,h,l,x)=\frac{\partial^2}{\partial h \partial l}\psi_{h,l,T}(x)=(1_{-l\le x\le h})\sum^{+\infty }_{n=-\infty }{4n(ng''_{T}(2n(h+l)-x)}\]
\[-\sum^{+\infty }_{n=-\infty }{(n-1)g''_{T}(2n(h+l)+x-2h))},\]

\par
\hfill \break
\par
\hfill \break
\noindent where 

\[g''_T(x)=\frac{x^2-T}{\sqrt{2\pi T^5}}e^{-\frac{x^2}{2T}}\]

\subsection{Density of the minimum, maximum, and terminal value for a Wiener process} For a Wiener process, a Brownian motion ${W_t, t>0}$, $W_0=0$, drift $\mu$, and volatility $\sigma$:

\[dW_{t}=\mu dt + \sigma dZ, \]
\[dZ \sim N(0,\sqrt{dt})\]

\hfill \break
\noindent the density of the minimum $l_T=min\{W_t\}_{t=0}^T$, maximum $h_T=max\{W_t\}_{t=0}^T$ and the terminal value $W_T$, can be presented as:

\[P(-l \le l_T \le h_T \le h; W_T\in [x,x+dx])=\psi^{\mu,\sigma}_{h,l,T}(x)dx\]
\[\psi^{\mu,\sigma}_{h,l,T}(x)=(1_{-l\le x\le h})e^{(\frac{\mu x}{\sigma^2} - \frac{\mu^2 T}{2\sigma^2})}\sum^{+\infty }_{n=-\infty }{(g_{\sigma^2 T}(2n(h+l)-x)-g_{\sigma^2 T}(2n(h+l)+x-2h))},\]

\[g_T(x)=\frac{1}{\sqrt{2\pi T  }}e^{-\frac{x^2}{2 T}}\]

\noindent Where we use \cite{rebholz}, similar formulation can be found in \cite{feller} p.430, but for a standard Brownian motion, without the drift and the volatility terms. To get to the density from this conditional formulation we derivate in respect to ${l}$ and ${h}$\footnote{The derivation was done in \cite{douady} however required some attention, please see the correction in 5.1}.

\[\psi''_{\mu,\sigma}(T,h,l,x)=\frac{\partial^2}{\partial h \partial l}\psi^{\mu,\sigma}_{h,l,T}(x)=(1_{-l\le x\le h})e^{(\frac{\mu x}{\sigma^2} - \frac{\mu^2 T}{2\sigma^2})} \sum^{+\infty }_{n=-\infty }{4n(ng''_{\sigma^2 T}(2n(h+l)-x)}\]
\[-\sum^{+\infty }_{n=-\infty }{(n-1)g''_{\sigma^2 T}(2n(h+l)+x-2h))},\]

\par
\hfill \break
\par
\hfill \break
\noindent where 

\[g''_T(x)=\frac{x^2-T}{\sqrt{2\pi T^5}}e^{-\frac{x^2}{2T}}\]

\subsection{Density of the range and terminal value for Wiener process} To find the join density of range and terminal value we integrate the trivariate density obtained earlier with respect with its minimum (or maximum):

\[\tilde{f}_{\mu,\sigma}(T,r,x)=\int^{r}_0\psi''_{\mu,\sigma}(T,h,r-h,x)dh\]

\[\psi''_{\mu,\sigma}(T,h,r-h,x)=(1_{h-r\le x\le h})e^{(\frac{\mu x}{\sigma^2} - \frac{\mu^2 T}{2\sigma^2})}\sum^{+\infty }_{n=-\infty }{4n^2g''_{\sigma^2 T}(2nr+x)-4n(n+1)g''_{\sigma^2 T}(2nr-x+2h)},\]

\[\psi''_{1,\mu,\sigma}(T,h,r-h,x) = 1_{h-r\le x\le h}e^{(\frac{\mu x}{\sigma^2} - \frac{\mu^2 T}{2\sigma^2})}\sum^{+\infty }_{n=-\infty }{4n^2g''_{\sigma^2 T}(2nr+x)},\]
\[\psi''_{2,\mu,\sigma}(T,h,r-h,x)=1_{h-r\le x\le h}e^{(\frac{\mu x}{\sigma^2} - \frac{\mu^2 T}{2\sigma^2})}\sum^{+\infty }_{n=-\infty }{-4n(n+1)g''_{\sigma^2 T}(2nr-x+2h)},\]
\[\psi''_{\mu,\sigma}(T,h,r-h,x)=\psi''_{1,\mu,\sigma}(T,h,r-h,x)+\psi''_{2,\mu,\sigma}(T,h,r-h,x),\]

\[\int^{r}_0\psi''_{1,\mu,\sigma}((T,h,r-h,x)dh = e^{(\frac{\mu x}{\sigma^2} - \frac{\mu^2 T}{2\sigma^2})}\int^{r}_0 1_{h-r\le x\le h}\sum^{+\infty }_{n=-\infty }{4n^2g''_{sigma^2 T}(2nr+x)dh},\]
\[= e^{(\frac{\mu x}{\sigma^2} - \frac{\mu^2 T}{2\sigma^2})}\int^{r}_0 1_{x\le h\le r+x}\sum^{+\infty }_{n=-\infty }{4n^2g''_{sigma^2 T}(2nr+x)dh},\]
\[= e^{(\frac{\mu x}{\sigma^2} - \frac{\mu^2 T}{2\sigma^2})}\sum^{+\infty }_{n=-\infty }{4n^2g''_{sigma^2 T}(2nr+x)}\int^{r}_0 1_{x\le h\le r+x}dh,\]
\[= e^{(\frac{\mu x}{\sigma^2} - \frac{\mu^2 T}{2\sigma^2})}\sum^{+\infty }_{n=-\infty }{4n^2g''_{sigma^2 T}(2nr+x)}(r-|x|),\]

\[\int^{r}_0\psi''_{2,\mu,\sigma}(T,h,r-h,x)dh = e^{(\frac{\mu x}{\sigma^2} - \frac{\mu^2 T}{2\sigma^2})}\int^{r}_0 1_{h-r\le x\le h}\sum^{+\infty }_{n=-\infty }{-4n(n+1)g''_{\sigma^2 T}(2nr-x+2h)},\]
\[= e^{(\frac{\mu x}{\sigma^2} - \frac{\mu^2 T}{2\sigma^2})}\int^{r}_0 1_{x\le h\le r+x}\sum^{+\infty }_{n=-\infty }{-4n(n+1)g''_{\sigma^2 T}(2nr-x+2h)},\]
\[= e^{(\frac{\mu x}{\sigma^2} - \frac{\mu^2 T}{2\sigma^2})}\sum^{+\infty }_{n=-\infty }\int^{r}_0 1_{x\le h\le r+x}{-4n(n+1)g''_{\sigma^2 T}(2nr-x+2h)}dh,\]
\[= e^{(\frac{\mu x}{\sigma^2} - \frac{\mu^2 T}{2\sigma^2})}\sum^{+\infty }_{n=-\infty }{-2n(n+1)(g'_{\sigma^2 T}(2nr+2r-|x|)-g'_{\sigma^2 T}(2nr+|x|))},\]

\noindent where
\[g'_T(x)=\frac{-x}{\sqrt{2\pi T^3}}e^{-\frac{x^2}{2T}}.\]

\par
\hfill \break
\par
\hfill \break
\[\tilde{f}_{\mu,\sigma}(T,r,x)=(1_{|x| \leq r})e^{(\frac{\mu x}{\sigma^2} - \frac{\mu^2 T}{2\sigma^2})}\sum^{+\infty }_{n=-\infty}{(4n^2g''_{\sigma^2 T}(2nr+x)(r-|x|)-}\]
\[ 2n(n+1)(g'_{\sigma^2 T}(2nr+2r-|x|)-g'_{\sigma^2 T}(2nr+|x|))\]

\begin{figure}[H]
  \includegraphics[width=\linewidth]{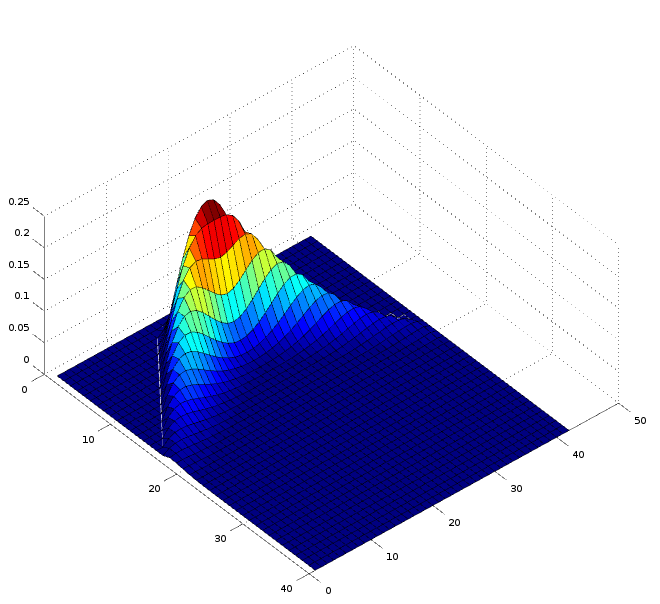}
  \caption{Join density of range and terminal value $\tilde{f}_{0,1}(T,r,x)$}
  \label{fig:boat1}
\end{figure}

\subsection{Density of range scaled terminal value of a standard Brownian Motion ($\mu=0,\sigma=1$)} To get to the range scaled terminal value we have the cumulative function as:

\[S(t,a) = P(\frac{x_{t}}{r_{t}}<a), a\in[-1,1]\] 
given the joint density function $\tilde{f}(t,r,x)$ of range $r$ and terminal value $x$ after time $t$, we can write:
\[S(t,a) = P(x_{t}<a{r_{t}})=\int_0^{\infty}{\int_{-r}^{ar}{\tilde{f}_{0,1}(t,r,x)}dx}dr\]
As we are interested in the density we can write:
\[s(t,a)=\frac{\partial S(t,a)}{\partial a} =\int_0^{\infty}{\frac{\partial}{\partial a}\bigg[{\tilde{F}_{0,1}(t,r,ar)-\tilde{F}_{0,1}(t,r,-r)}\bigg]}dr\]
\[=\int_0^{\infty}{\tilde{f}_{0,1}(t,r,ar)rdr}\]
where $\tilde{F}_{0,1}(t,r,x)$ is the primitive with respect to $x$.
\par
\[\tilde{f}_{0,1}(t,r,ar)=\sum^{+\infty }_{n=-\infty}{(4n^2g''_{T}(2nr+ar)(r-r|a|)-}\]
\[ 2n(n+1)(g'_{T}(2nr+2r-r|a|)-g'_{T}(2nr+r|a|))\]
for $a \in [-1,1],$
\[\tilde{f_1}(t,r,ar)=\sum^{+\infty }_{n=-\infty}{4n^2g''_T(2nr+ar)(r-r|a|)}\]
\[\tilde{f_2}(t,r,ar)=\sum^{+\infty }_{n=-\infty}{-2n(n+1)(g'_T(2nr+2r-r|a|)}\]
\[\tilde{f_3}(t,r,ar)=\sum^{+\infty }_{n=-\infty}{2n(n+1)(g'_T(2nr+r|a|)}\]
Integrating each of the factors $\tilde{f_1}(t,r,ar), \tilde{f_2}(t,r,ar), $ and $\tilde{f_3}(t,r,ar)$, we have to treat specially the cases when $n=0$, $n=-1$, as well as $a=0$:
\[\int_0^{\infty}\tilde{f_1}(t,r,ar)rdr=\sum^{+\infty }_{n=-\infty}\int_0^{\infty}{4n^2g''_T(2nr+ar)(1-|a|)r^2dr}\]
\[=\sum^{+\infty }_{n=-\infty}\frac{4n^2(1-|a|)}{((2H(-n)-1)a+2|n|)^3}\]

\[\int_0^{\infty}\tilde{f_2}(t,r,ar)rdr=\sum^{+\infty }_{n=-\infty}\int_0^{\infty}{-2(n+1)ng'_T(2nr+2r-r|a|)rdr}\]
\[=\sum^{+\infty }_{n=-\infty}\frac{(2H(n)-1)n(n+1)}{(-|a|+2n+2)^2}\]

\[\int_0^{\infty}\tilde{f_3}(t,r,ar)rdr=\sum^{+\infty }_{n=-\infty}\int_0^{\infty}{2(n+1)ng'_T(2nr+r|a|)rdr}\]
\[=\sum^{+\infty }_{n=-\infty}\frac{-(2H(n)-1)n(n+1)}{(|a|+2n)^2}\]
Finally  as the derivation not depend on $t$ we can write:
\[s(t,a)=\sum^{+\infty }_{n=-\infty}\bigg[\frac{4n^2(1-|a|)}{((2H(-n)-1)|a|+2|n|)^3}+ \frac{(2H(n)-1)n(n+1)}{(-|a|+2n+2)^2} - \frac{(2H(n)-1)n(n+1)}{(|a|+2n)^2}  \bigg]\]
\[s(a)=s(t,a)\]

\begin{figure}[H]
  \includegraphics[width=\linewidth]{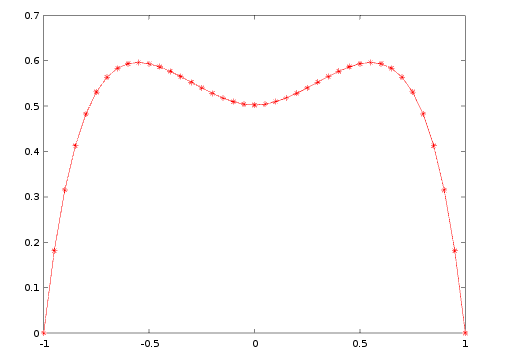}
  \caption{Density of range scaled terminal value of Brownian motion, given by $s(a)$}
  \label{fig:boat2}
\end{figure}
For $n=0$, the expression under the sum is equal to 0, even though it is not trivial to see in the case of the last summand under condition $a=0$, one can easily see going one step backward to the integral stage. From this point it is not even clear where the series converges. To simplify the formula we can eliminate the negative infinite limit by adding up the terms for negative $n$ with the positive $n$. After some transformations we can get to the following formulation of the range scale terminal value density:

\[s(a)=\sum^{+\infty }_{n=1}\bigg[\frac{16(1-|a|)n^3(4n^2+3a^2)}{(4n^2-a^2)^3}+ \frac{(8n^3|a|-8n^3-2na^2) }{(2n+|a|)^2(|a|-2n)^2} \bigg]\]
\[-\sum^{+\infty }_{n=2}\bigg[\frac{n(-8(n^2-1)|a|-2a^2+8(n^2-1))}{(-|a|-2n+2)^2(-|a|+2n+2)^2} \bigg]+\frac{2}{(4-|a|)^2}\]

\[s(a)=\sum^{+\infty }_{n=1}\bigg[\frac{16(1-|a|)n^3(4n^2+3a^2)}{(4n^2-a^2)^3}-\frac{(8n^3(1-|a|)+2na^2) }{(4n^2-a^2)^2} \bigg]\]
\[-\sum^{+\infty }_{n=2}\bigg[\frac{n(8(n^2-1)(1-|a|)-2a^2)}{(4n^2-(|a|-2)^2)^2} \bigg]+\frac{2}{(4-|a|)^2} \]

\par
\hfill \break
\par

\noindent To prove the convergence of the series we do the following transformation:

\[0 \leq s(a) \leq \sum^{+\infty }_{n=1}\bigg[\frac{16(1-|a|)n^3(4n^2+3a^2)}{(4n^2-a^2)^3}-\frac{(8n^3(1-|a|)+2na^2) }{(4n^2-a^2)^2} \bigg]\]
\[-\sum^{+\infty }_{n=2}\bigg[\frac{n(-8(n^2-1)(1-|a|)-2a^2)}{(4n^2-a^2)^2} \bigg]+\frac{2}{(4-|a|)^2} \leq \] 
\[\sum^{+\infty }_{n=1}\bigg[\frac{16(1-|a|)n^3(4n^2+3a^2)}{4n^2(4n^2-a^2)^2}-\frac{32n^5(1-|a|)+2na^2) }{4n^2(4n^2-a^2)^2} \bigg]\]
\[-\sum^{+\infty }_{n=2}\bigg[\frac{32n^3(n^2-1)(1-|a|)-2a^2)}{4n^2(4n^2-a^2)^2} \bigg]+\frac{2}{(4-|a|)^2} \leq \] 

\[\sum^{+\infty }_{n=2}\bigg[\frac{64(1-|a|)n^5+48n^3(1-|a|)a^2)}{4n^2(4n^2-a^2)^2}-\frac{32n^5(1-|a|)+2na^2) }{4n^2(4n^2-a^2)^2} \bigg]\]
\[-\sum^{+\infty }_{n=2}\bigg[\frac{32n^3(n^2-1)(1-|a|)-2a^2)}{4n^2(4n^2-a^2)^2} \bigg]+\frac{2}{(4-|a|)^2} + 1 = \] 

\[\sum^{+\infty }_{n=2}\bigg[\frac{80n^3(1-|a|)a^2-2na^2+2a^2}{4n^2(4n^2-a^2)^2}\bigg]+\frac{2}{(4-|a|)^2} + 1\] 

\noindent From here it is easy to see that the series converges.

\section{Empirical Results}
In this section we present the two studies, the first where we simulate paths from AR(1) process with different autocorrelation parameter and compare them with each other and in the second section we take market data and present the respective QQ plots.
\subsection{Simulating trend and mean reverting series} For the simulation of trend and mean reverting series we use AR(1) process as following:

\[W_t = \rho W_{t-1} + \epsilon_{t}\]

\noindent For $\rho=1$ this is Brownian motion examined previously, for $\rho<1$ the series is stationary while for trending and non-stationary behavior we can use $\rho>1$.

\begin{figure}[H]
\begin{minipage}{0.47\textwidth}
    \includegraphics[width=\linewidth]{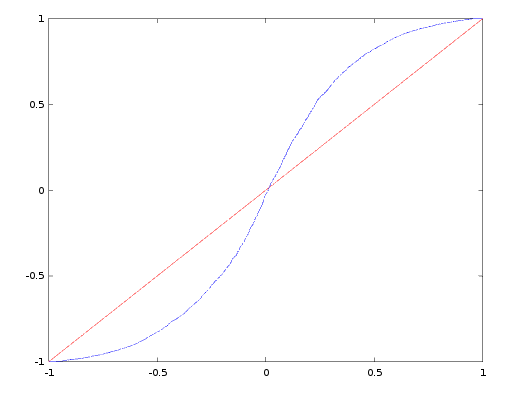}
    \end{minipage}
    \hspace{\fill} 
    \begin{minipage}{0.47\textwidth}
    \includegraphics[width=\linewidth]{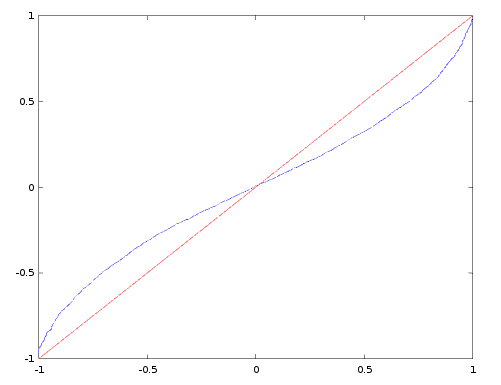}
    \end{minipage}

    \vspace*{1cm} 

\caption{QQ plots for tredning (left, $\rho=1.0005$) and mean reverting(right, $\rho=0.9995$) series against random walk (red line, $\rho=1$).} \label{fig:4pics}
\end{figure}
For trending series, the tails are more heavy while the middle is quite slim, while in mean revering series, the tails are flatter and most of the paths tend to end up more in the middle of its range. In other words, the bimodal humps of the s-density for mean revering series tend to become one hump. In the following figure we present the shapes of two mean revering series, the red one is simulated with paths of length $T=1500$ while the blue one is simulated with paths of length $T=500$ both with $\rho=0.9999$. The time invariance of the series with $\rho=1$ is not valid for $\rho \neq 1$.

 \begin{figure}[H]
  \includegraphics[width=\linewidth]{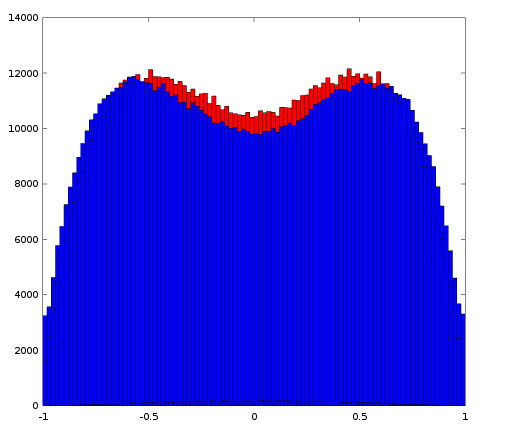}
  \caption{Simulation of AR(1) mean reverting process with different length of paths, in other words with different T.}
  \label{fig:boat3}
\end{figure}

\subsection{Market data} 
Here I consider two groups of currency pairs. In the first group currency pairs which are directly traded, pairs like EURUSD, USDJPY, EURGBP, GBPUSD, USDCHF, and USDCAD. The second group of currency pairs, are not directly traded,  commodity currency pair crosses such as AUDCAD, CADNOK, NZDNOK, NZDCAD, AUDNZD, CADNZD. To construct s-stat from empirical data we use the following:

\[s-stat=\frac{log(P_{close})-log(P_{open})}{log(P_{high})-log(P_{low})}\]

 \begin{figure}[H]
  \includegraphics[width=\linewidth]{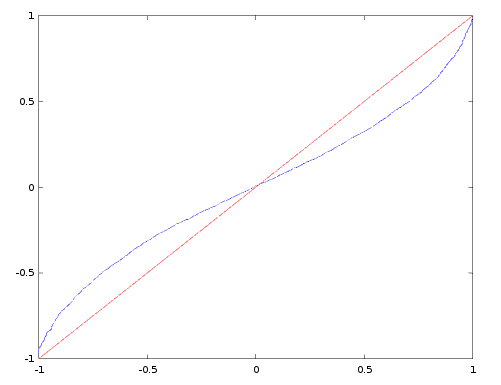}
  \caption{QQ Plot of group 1 (directly traded FX pairs) vs group 2 (indirectly traded FX pairs}
  \label{fig:boat3}
\end{figure}

\section{Conclusion}
The s-density obtained is bounded from -1 to 1 as expected and symmetrical around zero as expected. Another less trivial observation of this derivation and integration endeavour is that the range scaled terminal value is time invariant. It makes an intuitive sense, as the variance of the terminal value is linear in time as well as the expected range. This fact was quite expected and when the variable $t$ vanished going from $\tilde{f}(t,r,x)$ to $s(a)$ it was indicating that the derivation is on the right track. The other observation is less intuitive, the s-density is bimodal. An assumption which will require further investigation that in the limit for any $\rho<1$ the bimodality vanishes with T growing to infinity. In this case we can think of this bimodality as a bifurcation from stationary to non-stationary series.

\noindent

\section{Appendix}
\subsection{Correction} 
In the paper \cite{douady} the density of minimum, maximum, and the terminal value has already been derived, however there is a small correction need to be made to the formula (4.3). Here we need to correct for the fact that the range ${\delta}$ is not a difference between ${h}$ and ${l}$ but rather a sum, where the minimum value is smaller than ${-l}$ like in \cite{feller}. 

\[\psi_{h,l,T}(x)=1_{-l\le x\le h}\sum^{+\infty }_{n=-\infty }{(g_T(2n(h+l)-x)-g_T(2n(h+l)+x-2h))},\]

\noindent Given the symmetry of ${g_T}$, ${g''_T}$, as well as the sum itself (we can flip the sign of ${n}$, to go from ${n}$ to ${-n}$), we can get from Feller's formulation above to Douady's:

\[\psi_{h,l,T}(x)=1_{-l\le x\le h}\sum^{+\infty }_{n=-\infty }{(g_T(2(-n)(h+l)-x)-g_T(2(-n)(h+l)+x-2h))},\]
\[=1_{-l\le x\le h}\sum^{+\infty }_{n=-\infty }{(g_T(2n(h+l)+x)-g_T(2n(h+l)-x+2h))},\]

\noindent Finally,
\[\frac{\partial^2}{\partial h \partial l}\psi_{h,l,T}(x)=1_{-l\le x\le h} \sum^{+\infty }_{n=-\infty }{4n(ng''_T(2n(h+l)+x)}\]
\[-\sum^{+\infty }_{n=-\infty }{(n+1)g''_T(2n(h+l)-x+2h))},\]

\noindent Which is the exact formula in \cite{douady} except for the sign.

\subsection{Matlab code to generate s-density}
 
\begin{verbatim}
function y = s_density(a)
s = 0;

n = 1;
s = s + (1-abs(a))*16 * n^3 * (4 * n^2+3 * a^2) / ((4*n^2-a^2)^3);
s = s + (8*n^3*abs(a)-8*n^3-2*n*a^2) /(((2*n+abs(a))*(abs(a)-2*n))^2);
s = s + 2/(4-abs(a))^2;

for n = 2:100

s = s + (1-abs(a))*16 * n^3 * (4 * n^2+3 * a^2) / ((4*n^2-a^2)^3);
s = s + (8*n^3*abs(a)-8*n^3-2*n*a^2) /(((2*n+abs(a))*(abs(a)-2*n))^2);
s = s - n*(-8*(n^2-1)*abs(a)-2*a^2+8*(n^2-1))/(((-abs(a)-2*n+2)^2)*(((-abs(a)+2*n+2)^2)));
  
end
y = s

\end{verbatim}


\begin{thebibliography}{9}



\bibitem{trivariate}Choi, Byoung Seon, and Jeong Ho Roh. "On the trivariate joint distribution of Brownian motion and its maximum and minimum." Statistics \& Probability Letters 83.4 (2013): 1046-1053.

\bibitem{feller}Feller, William. "The asymptotic distribution of the range of sums of independent random variables." The Annals of Mathematical Statistics (1951): 427-432.

\bibitem{feller2}Feller, Willliam. An introduction to probability theory and its applications. Vol. 2. John Wiley \& Sons, 2008.

\bibitem{douady}Raphael, Douady. Closed form formulas for exotic options and their lifetime distributions. International Journal of Theoretical and Applied Finance, Vol.2, No. 1 (1998)

\bibitem{rebholz}He, H., Keirstead, W. P. and Rebholz, J. (1998), Double Lookbacks. Mathematical Finance, 8: 201–228. doi:10.1111/1467-9965.00053

\end{thebibliography}
\end{document}